\documentclass[a4paper,12pt]{article}

\usepackage{amsmath}
\usepackage{amsthm}
\usepackage{amssymb}  
\usepackage{latexsym} 
\usepackage[all]{xy}
\usepackage{comment}
\usepackage{rotating}
\usepackage{slashbox}
\usepackage{multirow}
\usepackage{rotating}

\newcommand{\Q}{{\mathbb Q}}

\newcommand{\Norm}{\operatorname{Norm}}

\newcommand{\Z}{{\mathbb Z}}

\newfont{\wncyr}{wncyr10 at 12pt}
\newfont{\wncyrten}{wncyr10 at 10pt}

\newenvironment{Proof}{\par\noindent{\sc Proof:}}%
                      {\hspace*{\fill}\nobreak$\Box$\par\medskip}
                       {\hspace*{\fill}\nobreak$\Box$\par\medskip}
\newenvironment{myitemize}
{\begin{itemize}
\setlength{\itemsep}{1pt}
\setlength{\parskip}{0pt}
\setlength{\parsep}{0pt}}
{\end{itemize}}

\newtheorem{Proposition}{Proposition}[section]
\newtheorem{Theorem}[Proposition]{Theorem}
\newtheorem{Lemma}[Proposition]{Lemma}
\newtheorem{Corollary}[Proposition]{Corollary}
\newtheorem{Conjecture}[Proposition]{Conjecture}

\theoremstyle{definition}

\newtheorem{Example}[Proposition]{Example}

\newcounter{nootje}
\begin{document}
\normalsize
\title{On the Diophantine equation $\displaystyle cy^l=\frac{x^p-1}{x-1}$}
\author{Mohammad Sadek\\ Department of Mathematics and Actuarial Sciences\\ American University in Cairo\\ mmsadek@aucegypt.edu}
\date{}
\maketitle
\begin{abstract}{\footnotesize Let $\displaystyle p,\,c$ be distinct odd primes, and $l\ge 2$ an integer. We find sufficient conditions for the Diophantine equation \[cy^{l}=\Phi_p(x)=\frac{x^p-1}{x-1}=x^{p-1}+x^{p-2}+\ldots+1\] not to have integer solutions.}
\end{abstract}

\section{Introduction}

The solutions of the Nagell-Ljunggren equation $\displaystyle y^q=\frac{x^n-1}{x-1}$, where $q,n\ge 2$ are integers, have been the source for many conjectures. One of these is the following:
\begin{Conjecture}
The only solutions to the Diophantine equation $\displaystyle y^q=\frac{x^n-1}{x-1}$ in integers $x,y>1,n>2,q\ge 2$ are given by \[\frac{3^5-1}{3-1}=11^2,\;\frac{7^4-1}{7-1}=20^2,\textrm{ and }\frac{18^3-1}{18-1}=7^3.\]
\end{Conjecture}
The above conjecture has been solved completely for $q=2$. Furthermore, it has been proved if one of the following assumptions holds:
 \[3\mid n,\textrm{ or }4\mid n,\textrm{ or }q=3 \textrm{ and }n\not\equiv5\textrm{ mod }6.\] We moreover know that the Nagell-Ljunggren equation has no solutions with $x$ square. The main tools used to attack this Diophantine equation are effective Diophantine approximation, linear forms in $p$-adic logarithms, and Cyclotomic fields theory. For these results and more see \cite{BegeaudMignotte}, \cite{Mihailescu2} and \cite{Mihailescu1}.

In \cite{Bugeaud} the Diophantine equation $\displaystyle y^l=c\frac{x^n-1}{x-1}$ has been treated. A complete list of such Diophantine equations with integer solutions has been given, under the condition that $1\le c\le x\le 100$. A more general equation $a\frac{x^n-1}{x-1} = cy^l$ where $ac> 1$ has been considered in \cite{Shorey}. Our interest in the latter equation is when $a=1$.

In this note we will be concerned with the Diophantine equation $\displaystyle cy^l=\frac{x^p-1}{x-1}$, where $c,p$ are distinct odd primes and $l\ge 2$. We exhibit the existence of an infinite set of triples $(p,c,l)$ for which the mentioned Diophantine equation has no integer solutions. For example, this infinite set contains the set of triples $(p,c,l)$ where the Legendre symbol $\displaystyle\left(\frac{c}{p}\right)=-1$ and $l$ is even.

 The key idea is exploiting the following identity satisfied by the cyclotomic polynomial $\displaystyle\Phi_p(x)=\frac{x^p-1}{x-1}$  \[4\Phi_p(x)=A_p(x)^2-(-1)^{(p-1)/2}pB_p(x)^2,\] where $A_p(x),B_p(x)\in\Z[x].$ This identity goes back to Gauss, nevertheless the formulae describing $A_p(x)$ and $B_p(x)$ were given recently in \cite{Brent}. Using this identity we show that the existence of an integer solution to the equation in question implies the existence of a proper integer solution to some auxiliary Diophantine equation.

\section{Factorization of cyclotomic polynomials}

For an odd square-free integer $n>1$, and $|x|\le 1$ define \[f_n(x)=\sum_{j=1}^{\infty}\left(\frac{j}{n}\right)\frac{x^j}{j},\]
where $\displaystyle\left(\frac{j}{n}\right)$ is the Jacobi symbol of $j$ mod $n$. We state Theorem 1 of \cite{Brent}.

\begin{Theorem}
\label{thm1}
Let $n>3$ be an odd square-free integer. Consider the Gauss's identity $\displaystyle 4\Phi_n(x)=A_n(x)^2-(-1)^{(n-1)/2}nB_n(x)^2$, where $A_n(x),B_n(x)\in\Z[x]$. If $n\equiv 1$ mod $4$, then
\begin{eqnarray}
A_n(x)=2\sqrt{\Phi_n(x)}\cosh\left(\frac{\sqrt{n}}{2}f_n(x)\right),\;B_n(x)=2\sqrt{\frac{\Phi_n(x)}{n}}\sinh\left(\frac{\sqrt{n}}{2}f_n(x)\right).\nonumber
\end{eqnarray}
 If $n\equiv 3$ mod $4$, then
\begin{eqnarray}
A_n(x)=2\sqrt{\Phi_n(x)}\cos\left(\frac{\sqrt{n}}{2}f_n(x)\right),\;B_n(x)=2\sqrt{\frac{\Phi_n(x)}{n}}\sin\left(\frac{\sqrt{n}}{2}f_n(x)\right).\nonumber
\end{eqnarray}
\end{Theorem}

\section{An auxiliary Diophantine equation}
\label{sec2}

The results of this section are motivated by Proposition 8.1 of \cite{DarmonGran}.

By a proper solution $(x_0,y_0,z_0)$ to the Diophantine equation $\displaystyle ax^p+by^q=cz^r$, we mean three integers $x_0,y_0,z_0$ such that $\displaystyle ax_0^p+by_0^q=cz_0^r$ and $\gcd(x_0,y_0,z_0)=1$.

We state the following result on local solutions to $cy^l=x^2\pm p z^2$ where $c,p$ are distinct odd primes and $l\ge 2$.

\begin{Proposition}
\label{prop:local}
 There are proper local solutions to \[\alpha ^2 cy^l=x^2 \pm p z^2,\;\alpha\in\{1,2\},\] at every prime if and only if the Legendre symbol $\left(\displaystyle\frac{\mp p}{c}\right)=1$; and, when $l$ is even we have $\displaystyle \left(\frac{c}{p}\right) = 1$.
\end{Proposition}
\begin{Proof}
The given conditions are clearly necessary. Now we need to prove they are sufficient. We use the fact that if $q\nmid 2cp$, then there are $q$-adic
 integer solutions to $x^2\pm pz^2=\alpha^2 c$, so take $(x,1,z)$. For the prime $c$, since $\left(\displaystyle\frac{\mp p}{c}\right)=1$, there are $c$-adic integer solutions to $x^2=\mp p$, so take $(x,0,1)$. For the prime $p$, if $l$ is odd, take $(\alpha c^{(l+1)/2},c,0)$; if $l$ is even, hence $\left(\displaystyle\frac{c}{p}\right)=1$, then there is a $p$-adic integer satisfying $x^2=\alpha^2 c$, and we take $(x,1,0)$. For the prime $2$, the equation becomes $x^2- z^2=\alpha^2y^l$, so we can lift the solution $(1,0,1)$ mod $2$ to a $2$-adic integer solution.
\end{Proof}

\begin{Proposition}
\label{prop:Darmon}
Let $\displaystyle p,c$ be distinct odd primes, and $l\ge 2$ be an integer. Set $\displaystyle\delta=(-1)^{(p-1)/2}$. If the Diophantine equation \[\displaystyle \alpha^2 cy^l=x^2-\delta pz^2,\;\alpha\in\{1,2\},\] has a proper solution with $y$ being odd and $\gcd(x,y)=1$, then there exist coprime ideals $I, J$ in $\displaystyle\Q(\sqrt{\delta p})$ with $I J = (\alpha^2 c)$, whose ideal classes are $l$-th powers inside the class group of $\displaystyle\Q(\sqrt{\delta p})$.
\end{Proposition}
\begin{Proof}
Suppose $(x,y,z)$ is a proper solution to $\displaystyle \alpha^2cy^l=x^2-\delta pz^2$ where $y$ is odd and $\gcd(x,z)=1$. Now considering the latter as ideal equation, we have
 \[(\alpha^2c)(y)^l=(x-\sqrt{\delta p}\;z)(x+\sqrt{\delta p}\;z).\]
 Now the ideal $\mathfrak{a}=(\displaystyle x-\sqrt{\delta p}\;z,x+\sqrt{\delta p}\;z) \mid (2x,2\sqrt{\delta p}, \alpha^2cy^l)=(2,\alpha)$.
 \begin{myitemize}
 \item[1)] If $\alpha=1$, then \[(x-\sqrt{\delta p}\;z)=I L_1^l,\;(x+\sqrt{\delta p}\;z)=J L_2^l,\]
 where $I J=(c)$ and $L_1L_2=(y)$. This implies that the ideal classes of $I$ and $J$ are both $l$-th powers inside the class group of
 $\displaystyle \Q(\sqrt{\delta p})$.
 \item[2)] If $\alpha=2$, then both $x,z$ are odd. This will yield a contradiction when $p\equiv \pm1$ mod $8$. This follows from the fact that $4cy^l=x^2-\delta pz^2\equiv 0$ mod $8$ when $p\equiv \pm 1$ mod $8$.

     When $p\equiv \pm 5$ mod $8$, the ideal $(2)$ is prime inside $\displaystyle \Q(\sqrt{\delta p})$ because $\delta p\equiv 5$ mod $8$. If $\mathfrak{a}=(2)$, then $\displaystyle 2\mid (x-\sqrt{\delta p}z)$ which implies that $2\mid x,z$, a contradiction. Thus $\mathfrak{a}=1$, and we argue like in the first case.
\end{myitemize}
\end{Proof}

\section{The equation $\displaystyle cy^l=\frac{x^p-1}{x-1}$}

We start by stating the following elementary lemma.

\begin{Lemma}
\label{lem:gcd}
Let $a\in\Z$ and $p$ be an odd prime.
\begin{myitemize}
\item[i)] $\Phi_p(a)$ is odd.
\item[ii)] Set $d=\gcd(A_p(a),B_p(a))$. Then $d\in\{1,2\}$. If $p\equiv \pm 1$ mod $8$, then $d=2$.
\end{myitemize}
\end{Lemma}
\begin{Proof}
i) Since $\Phi_p(a)\equiv 1$ mod $a$, hence if $a$ is even, $\Phi_p(a)$ is odd. If a is odd,
then $\Phi_p(a)\equiv \Phi(1)=p$ mod $2$.

ii) Assume that $q\mid d$, where $q> 1$ is an odd prime. We will write $\tilde{a}$ for the reduction of $a$ mod $q$.

 If $q\ne p$, then $\displaystyle (x-\tilde{a})\mid A_p(x),\;B_p(x)$ mod $q$ because $q\mid A_p(a)$ and $B_p(a)$. Hence $\displaystyle (x-\tilde{a})^2\mid \Phi_p(x)$ mod $q$. The latter statement contradicts the fact that $x^p-1$ has no multiple factors mod $q$ when $\gcd(q,p)=1$.

 If $q=p$, then $p^2\mid \Phi_p(a)$. In particular $\tilde{a}^p\equiv 1$ mod $p$. Fermat's Little Theorem yields that there is a $\lambda\in\Z$ such that $a= 1+\lambda p.$ So
 \begin{eqnarray*}
 \Phi_p(a)&=&\sum_{i=0}^{p-1} a^{i}=\sum_{i=0}^{p-1}( 1+\lambda p)^i\\
 &\equiv& p +\lambda p\sum_{i=0}^{p-1}i\equiv p\textrm{ mod }p^2,
 \end{eqnarray*}
 which contradicts that $p^2\mid\Phi_p(a)$. We conclude that $d\mid 2$.

 Now we assume $p\equiv\pm 1$ mod $8$. Assume on the contrary that $2\nmid d$. This implies that both $A_p(a)$ and $B_p(a)$ are odd as $4\mid A_p^2(a)-(-1)^{(p-1)/2}pB_p^2(a)$. A direct calculation shows that if $A_p(a),B_p(a)$ are both odd, then \[4\Phi_p(a)=A_p^2(a)-(-1)^{(p-1)/2}pB_p^2(a)\equiv1-(-1)^{(p-1)/2}p\equiv 0 \textrm{ mod }8,\] which contradicts (i).
\end{Proof}

\begin{Corollary}
\label{cor1}
Let $p,\,c$ be distinct odd primes. Let $l\ge2$ be an integer. Assume that $(a,b)$ is an integer solution to the Diophantine equation $\displaystyle cy^l=\Phi_p(x)$. Then there exists an integer solution $(x,y,z)$, where $\gcd(x,z)=1$ and $y$ is odd, to a Diophantine equation of the form \[\alpha^2cy^l=x^2-(-1)^{(p-1)/2}pz^2,\;\alpha\in\{1,2\}.\]
In the case $p\equiv \pm1$ mod $8$, one has $\alpha=1$.
\end{Corollary}
\begin{Proof}
One has $4 c\,b^l=4\Phi_p(a)=A_p(a)^2-(-1)^{(p-1)/2}pB_p(a)^2$, where $A_p(x),B_p(x)\in\Z[x]$ and $d=\gcd(A_p(a),B_p(a))\mid 2$, Lemma \ref{lem:gcd}. If $d=1$, then $(A_p(a),b,B_p(a))$ is a proper solution to $\displaystyle 4cy^l=x^2-(-1)^{(p-1)/2}z^2$. If $d=2$, then $(A_p(a)/2,b,B_p(a)/2)$ is a proper solution to $\displaystyle cy^l=x^2-(-1)^{(p-1)/2}z^2$. Observe that if $p\equiv\pm1$ mod $8$, then $d=2$, Lemma \ref{lem:gcd} (ii).
\end{Proof}

Now we state our main result which says that there is an infinite number of triples $(c,p,l)$ such that $\displaystyle cy^l=\Phi_p(x)$ has no integer solution.

\begin{Theorem}
\label{thm:main}
Let $\displaystyle p,\,c$ be distinct odd primes, and $l\ge 2$ an integer. Set $\delta=(-1)^{(p-1)/2}$. If the triple $\displaystyle (p,c,l)$ satisfies one of the following conditions:
\begin{myitemize}
\item[i)] $\displaystyle \left(\frac{\delta p}{c}\right)=-1$;
\item[ii)] $\displaystyle \left(\frac{c}{p}\right)=-1$, and $l$ is even;
\item[iii)] There exist no ideals $I,J$ whose ideal classes are $l$-th powers in the class group of $\displaystyle\Q(\sqrt{\delta p})$ and satisfy $(\alpha^2c)=IJ$, where
    \[\alpha\in\left\{\begin{array}{lr}\{1\}&\textrm{ if }p\equiv \pm1\textrm{ mod }8\\\{1,2\} & \textrm{ if }p\equiv \pm3 \textrm{ mod }8\end{array}\right.\]
\end{myitemize}
then the Diophantine equation \[cy^l=\Phi_p(x)=\frac{x^p-1}{x-1}=x^{p-1}+x^{p-2}+\ldots+x+1\] has no integer solutions.
\end{Theorem}
\begin{Proof}
Assume on the contrary that there exists a proper integer solution to $c y^l=\Phi_p(x)$. This implies the existence of a proper integer solution to $\alpha^2cy^l=x^2-\delta p z^2$, see Corollary \ref{cor1}. Hence we have a contradiction in (i) and (ii), see Proposition \ref{prop:local}. Furthermore one has a contradiction in case (iii) obtained using Proposition \ref{prop:Darmon}.
\end{Proof}

Parts (i) and (ii) of the above theorem provide an infinite family of Diophantine equations with no integer solutions. For example
\[13y^{2l}=\Phi_{137}(x)=x^{136}+\ldots+x+1\] has no integer solutions because $\displaystyle \left(\frac{13}{137}\right)=-1$.

In the following example we show that (iii) of Theorem \ref{thm:main} can be used to find explicit triples $(c,l,p)$ such that the Diophantine equation $cy^l=\Phi_p(x)$ has no integer solutions.

\begin{Example}
The Diophantine equation\[3y^{5k}=\Phi_{47}(x)=x^{46}+x^{45}+\ldots+x+1,\;k\ge1,\] has no integer solutions. We have $47\equiv-1$ mod $8$ and $(3)=\mathfrak{p}\mathfrak{p}'$ in the ring of integers of $\displaystyle \Q(\sqrt{-47})$. The class number of $\displaystyle \Q(\sqrt{-47})$ is $5$. The ideal class $[\mathfrak{p}]$ of $\mathfrak{p}$ can not be a fifth power inside the ideal class group of $\displaystyle \Q(\sqrt{-47})$ because $[\mathfrak{p}]$ generates the ideal class group.
\end{Example}

\bibliographystyle{plain}
\footnotesize
\bibliography{thesisreferences}
\end{document}